\documentclass[12pt]{amsart}
\usepackage{amsfonts}
\usepackage{amsmath}
\usepackage{amssymb}
\usepackage{amsthm}

\input xy
\xyoption {all}

\def \m {{\mathfrak M}}
\def \ms{{\mathbf M}}
\def \ts{{\mathbf \Theta}}

\newcommand{\comment}[1]{}
\newtheorem*{thma}{Theorem 1A}
\newtheorem {question}{Question}

\newtheorem{theorem}{Theorem}
\newtheorem {lemma}{Lemma}
\newtheorem{conjecture}{Conjecture}
\newtheorem {corollary}{Corollary}
\newtheorem {proposition}{Proposition}
\newtheorem {claim}{Claim}

\theoremstyle{definition}
\newtheorem{remark}{Remark}

\theoremstyle {definition}

\begin{document}
\baselineskip=16pt

\title[Generic strange duality for $K3$ surfaces]{Generic strange duality for $K3$ surfaces}
\author {Alina Marian}
\address {Department of Mathematics}
\address {University of Illinois at Chicago}
\email {alina@math.uic.edu}
\author {Dragos Oprea}
\address {Department of Mathematics}
\address {University of California, San Diego}
\email {doprea@math.ucsd.edu}
\address {Department of Mathematics}
\address{Faculty of Science, Kobe University, Kobe, 657, Japan}
\email {yoshioka@math.kobe-u.ac.jp}
\date{}

\begin{abstract}
Strange duality is shown to hold over generic $K3$ surfaces in a large number of cases. The isomorphism for elliptic $K3$ surfaces is established first via Fourier-Mukai techniques. Applications to Brill-Noether theory for sheaves on $K3$s are also obtained. The appendix written by Kota Yoshioka discusses the behavior of the moduli spaces under change of polarization, as needed in the argument.
\end {abstract}
\maketitle

\section{Introduction}

\subsection {The strange duality morphism} 

We consider moduli spaces of sheaves over $K3$ surfaces, and the strange duality map on spaces of generalized theta functions associated to them. 

To start, we recall the general geometric setting for strange duality phenomena. Let $(X, H)$ be a  smooth polarized complex projective surface. To give our exposition a uniform character, we assume that $X$ is simply connected. Let $v$ be a class in the topological $K$-theory $K_{\text{top}} (X)$ of the surface, and denote by $\m_v$ the moduli space of Gieseker $H$-semistable sheaves on $X$ of topological type $v$. 

The moduli space $\m_v$ carries natural line bundles which we now discuss. Consider the bilinear form on $K_{\text{top}} (X)$ given by
\begin{equation}
( v, w ) = \chi (v \cdot w), \, \, \text{for} \, \, v, w \in K_{\text{top}} (X),
\label{euler}
\end{equation}
where the  product in $K$-theory is used. Let $$v^{\perp} \subset K_{\text{top}} (X)$$ contain the $K$-classes orthogonal to $v$ relative to this form. When $\m_v$ consists of stable sheaves only\footnote{The homomorphism is defined in all generality from a more restricted domain.}, there is a group homomorphism
$$\Theta: v^{\perp} \to \, \text{Pic} \, \m_v, \, \, \, w \mapsto \Theta_w,$$
studied among others in \cite{lepotier}, \cite{jun}. If $\m_v$ carries a universal sheaf $$\mathcal E \to \m_v \times X,$$ and $w$ is the class of a vector bundle $F$, we have  
$$\Theta_w = \det {\mathbf R}p_! (\mathcal E \otimes q^{\star} F)^{-1}.$$
Here $p$ and $q$ are the two projection maps from $\m_v \times X.$ The theta line bundle is also defined in the absence of a universal sheaf, by
descent from the Quot scheme. 

We consider now two classes $v$ and $w$ in $K_{\text{top}} (X)$ satisfying
$$( v, w ) = 0.$$  
If the conditions 
$$H^2 (E \otimes F) = 0 ,\,\,\, \text{Tor}^{1}(E,F)=\text{Tor}^2(E,F)=0$$ hold in $\m_v \times \m_w$ away from codimension $2$, the locus 
\begin{equation}\label{thetaeq}\Theta = \{(E, F) \in \m_v \times \m_w \, \, \text{such that} \, \, H^0 (E \otimes F) \neq 0 \}\end{equation} should correspond to a divisor.
Furthermore, 
$${\mathcal O} (\Theta) = \Theta_w \boxtimes \Theta_v,$$ so $\Theta$ induces a map 
$${\mathsf D}: H^0 (\m_v, \Theta_w)^{\vee} \to H^0 (\m_w, \Theta_v).$$ 
The main query concerning this map is
\begin{question}
When nonzero, is $\mathsf D$ an isomorphism? 
\end{question}
While the question is  too naive  for an affirmative answer to be expected in this generality, the isomorphism was shown to hold for many pairs $(\m_v, \m_w)$ of moduli spaces of sheaves over either $K3$ or rational surfaces, cf.
\cite {A} \cite {D1} \cite {D2} \cite {G} \cite {OG2} \cite {S}. In all examples however,  one of the moduli spaces involved has small dimension and the other consists of rank $2$ sheaves.  A survey of some of the known results is contained in \cite {mo}. In this paper, we establish the isomorphism on moduli spaces over generic $K3$ surfaces for a large class of topological types of the sheaves, allowing in particular for arbitrarily high ranks and dimensions. The precise statements are as follows.

\subsection{Results} Let $(X, H)$ be a polarized $K3$ surface. We use as customary the Mukai vector $$v(E) = \text{ch} E \sqrt{\text{Todd}\, X} \,  \in \, H^{\star} (X, {\mathbb Z})$$ to express the topological type of a sheaf $E$ on $X$. 
We write $$v = v_0 + v_2 + v_4$$ to distinguish cohomological degrees in $v$, and set $$v^{\vee} = v_0 - v_2 + v_4.$$
Note also the Mukai pairing on cohomology: $$\langle v, w \rangle = \int_S v_2 w_2 - v_0 w_4 - v_4 w_0.$$  In terms of the pairing \eqref{euler}, we have 
$$(v, w) = \langle v, w^{\vee} \rangle = \langle v^{\vee}, w \rangle.$$
We assume that the moduli space $\m_{v}$ of Gieseker $H$-semistable sheaves of fixed Mukai vector $v$ consists only of stable sheaves. In this case,  $\m_v$ is an irreducible holomorphic symplectic manifold whose dimension is simply expressed in terms of the Mukai self-pairing of $v$, $$\dim \m_v = \langle v, v\rangle + 2 .$$

We show

\begin {theorem}\label{genericarb}
Assume $(X, H)$ is a generic polarized $K3$ surface with $\text{Pic}\, X = {\mathbb Z}\, H$, and consider orthogonal Mukai vectors $v$ and $w$ of ranks $r\geq 2$ and $s\geq 3$ such that
\begin {itemize}
\item [(i)] $c_1(v)=c_1(w)=H,$
\item [(ii)] $\chi (v) \leq 0$, $\chi(w)\leq0,$
\item [(iii)] $\langle v, v\rangle \geq 2 (r-1)(r^2+1), \,\,\langle w, w\rangle \geq 2 (s-1)(s^2+1).$
\end {itemize}
Then $${\mathsf D}:H^0(\m_v, \Theta_w)^{\vee}\to H^0(\m_w, \Theta_v)$$ is an isomorphism. 
\end {theorem}

The genericity means that the statement holds on a nonempty open subscheme of the moduli space of polarized $K3$s. We expect the result to be true for all $K3$s and we will pursue this aspect in future work. 

Similarly, in rank $2$, we prove 
\begin {thma} \label{generict} Assume $(X, H)$ is a generic polarized $K3$ surface of degree at least $8$, and consider orthogonal Mukai vectors $v$ and $w$ of rank $2$ such that 
\begin {itemize}
\item [(i)] $c_1(v)=c_1(w)=H,$
\item [(ii)] $\chi (v) \leq 0$, $\chi(w)\leq0.$
\end {itemize}
Then $${\mathsf D}:H^0(\m_v, \Theta_w)^{\vee}\to H^0(\m_w, \Theta_v)$$ is an isomorphism.
\end {thma}

The statement is obtained by degeneration to moduli spaces over a smooth elliptic $K3$ surface $X$ with a section. Better results are in fact available here. Let us assume that the fibers have at worst nodal singularities and that the N\'{e}ron-Severi group is 
$$\text{NS}\, (X) = {\mathbb Z} \sigma + {\mathbb Z} f,$$ where $\sigma$ and $f$ are the classes of the section and of the fiber respectively. For fixed Mukai vectors, we consider stability with respect to polarizations $H=\sigma+mf$ suitable in the sense of \cite {F}. We show

\begin{theorem}
\label{sd}
Let $v$ and $w$ be orthogonal Mukai vectors corresponding to sheaves of ranks $r$ and $s$ on $X$ with $r, s \geq 2$. Assume further that
\begin{itemize}\item [(i)]$c_{1}(v)\cdot f=c_{1}(w)\cdot f=1$, \item [(ii)] $\langle v, v\rangle + \langle w, w\rangle \geq 2(r+s)^{2}.$\end{itemize} Then the duality map $$\mathsf {D}:H^{0}(\m_{v}, \Theta_{w})^{\vee}\to H^{0}(\m_{w}, \Theta_{v})$$ is an isomorphism. 
\end{theorem}

Along the way we establish the following Brill-Noether result for sheaves on $K3$ elliptic surfaces.

\begin {theorem} \label{t3}Under the assumptions of Theorem \ref{sd}, the locus $\Theta$ has codimension $1$ in the product of moduli spaces $\m_v\times \m_w$. In particular, for a generic sheaf $E\in \m_v$, $$\Theta_{E}=\{F\in \m_w: h^0(E\otimes F)\neq 0\}$$ is a divisor in $\m_w$. 
\end {theorem}

The proofs use the fact that the moduli spaces $\m_v$ and $\m_w$ are birational to Hilbert schemes of points on $X$, 
\begin{equation}
\label{crucial}
\m_v \dashrightarrow X^{[a]}, \, \, \, \, \m_w \dashrightarrow X^{[b]}, \, \, \, \text{with} \, \, \, a = \frac {\langle v, v\rangle}{2} + 1, \, \, \, b =\frac {\langle w, w\rangle}{2} + 1. 
\end{equation} The birational maps \eqref{crucial} were described in \cite{ogrady} and were shown to be regular away from codimension $2$. Theorem \ref{sd} is then a consequence of the explicit identification of the theta divisor \eqref{thetaeq} with a divisor  in the product $X^{[a]} \times X^{[b]}$ known to induce strange duality. Specifically, for any line bundle $L$ on  $X$ with $\chi(L) = a+b$ and no higher cohomology, one can consider the divisor associated to the locus $$\theta_{L, a,b}=\{(I_Z, I_W)\, \, \text{such that} \, \, h^0 (I_Z \otimes I_W \otimes L) \neq 0 \} \subset X^{[a]} \times X^{[b]}.$$ 
Furthermore, observe the involution on the elliptic surface $X$, given by fiberwise reflection across the origin of  the fiber: 
$$p \in f \, \mapsto -p \in f.$$ The involution is defined away from the codimension 2 locus of singular points of fibers of $X$. It induces an involution on any Hilbert scheme of points on $X$,
defined outside a codimension 2 locus, $$\iota: X^{[a]} \dashrightarrow X^{[a]}, \, \, \, Z \mapsto \widetilde{Z}.$$
Consider the pullback $$\tilde{\theta}_{L,a,b} = (\iota\times 1)^{\star} {\theta}_{L,a,b},
$$ under the birational map
$$
\iota\times 1: \, X^{[a]} \times X^{[b]} \dashrightarrow X^{[a]} \times X^{[b]}, \, \, \, (Z, W) \, \mapsto (\widetilde Z, {W}).
$$ It can be shown that  
$$(\iota\times 1)^{\star} \theta_{L, a, b} = (1\times \iota)^{\star} \theta_{L, a, b}.$$
Viewing $\Theta$ as a locus in $X^{[a]} \times X^{[b]}$ via \eqref{crucial}, we prove 

\begin{theorem} \label{theta} There exists a line bundle $L$ on $X$, such that 
$${\Theta} = \tilde{\theta}_{L, a, b} \, \, \, \text{in the product} \, \, \,  X^{[a]}\times X^{[b]}.$$
\end{theorem}

$\theta_{L,a,b}$ is known to give an isomorphism on the associated spaces of sections on $X^{[a]}$ and $X^{[b]}$, cf. \cite{mo}. Therefore so does $\tilde{\theta}_{L,a,b},$ yielding Theorem \ref{sd}.

The identification of the two theta divisors of Theorem \ref{theta} is difficult even though the O'Grady birational isomorphism with the Hilbert scheme is explicit. To achieve it, we interpret the O'Grady construction by means of Fourier-Mukai transforms. We show that the Fourier-Mukai transforms of {\it generic} O'Grady sheaves are two-term complexes in the derived category, derived dual to ideal sheaves. More importantly, a careful analysis is necessary to keep track of the special loci where the generic description may fail.  

The same method gives results for arbitrary simply connected elliptic surfaces $$\pi:X\to \mathbb P^1$$ with a section and at worst nodal fibers. The dimension of the two complementary moduli spaces will be taken large enough compared to the constant $$\Delta= \chi(X, \mathcal O_X) \cdot \left((r+s)^2+(r+s)+2\right)-2(r+s).$$ We continue to assume that the polarization is suitable. We prove

\begin {theorem}\label{ellsur} Assume $v$ and $w$ are two orthogonal topological types of rank $r, s\geq 2$, such that \begin {itemize} \item [(i)] $c_1(v)\cdot f=c_1(w)\cdot f=1,$
\item [(ii)] $\dim \m_v+\dim \m_w\geq \Delta.$ \end {itemize} Then, $\Theta$ is a divisor in $\m_v\times \m_w$. 
\end {theorem}

We also propose

\begin {conjecture} \label{ellsurconj} Under the assumptions of Theorem \ref{ellsur}, $$\mathsf D:H^0(\m_v, \Theta_w)^{\vee}\to H^0(\m_w, \Theta_v)$$ is an isomorphism. 
\end {conjecture}

The conjecture is in fact established up to the statement that, in this new setting, the birational isomorphism of $\m_v$ and $\m_w$ with Hilbert schemes of points holds {\it away from codimension 2}. We believe this to be true. 

The paper is structured as follows. The main part of the argument concerns the case of elliptic $K3$ surfaces and is presented in Section $2$. The last part of Section 2 treats the case of arbitrary simply connected elliptic surfaces. Section $3$ explains generic strange duality via a deformation argument. The appendix written by Kota Yoshioka contains a discussion of change of polarization for higher-rank moduli spaces of sheaves over $K3$s.

\subsection {Acknowledgments} We thank Kieran O'Grady, Eyal Markman,  Michael Thaddeus, and Kota Yoshioka for helpful conversations and correspondence.  In particular, Kota Yoshioka suggested a simplification of the main calculation of Section 2.5. We are grateful to the NSF for financial support. 

\section{The theta isomorphism for elliptic $K3$ surfaces}

\subsection{O'Grady's construction} 
Keeping the notations of the introduction, we let $\pi: X \to {\mathbb P}^1$ be an elliptic $K3$ surface with a section $\sigma$, whose fibers have at worst nodal singularities. We have
$$\sigma^2 = -2, \, \, f^2 = 0, \, \, \sigma\cdot f  = 1.$$ 

We are concerned with sheaves on $X$ with Mukai vector of type 
$$ v = r + (\sigma + k f) + p\,\omega,$$ for some $k, p \in {\mathbb Z},$ with $\omega$ being the class of a point in $X$. We consider a $v$-{\it{suitable}} polarization $$H=\sigma+mf \, \, \, \text{for} \, \, \, m>>0.$$ This means that $H$ lies in a $v$-chamber of the ample cone of $X$ adjacent to the class $f$ of the fiber \cite{ogrady}. The moduli space $\m_v$ of $H$-semistable sheaves consists of slope-stable sheaves only, and the choice of $H$ ensures that $E \in \m_v$ is stable if and only if its restriction to a generic fiber is stable. 
The restriction to special fibers may be unstable, as we will show in Lemma \ref{spl1} below. 

As explained in \cite {ogrady}, we can inductively build the moduli spaces $\m_v$ as follows. Note first that tensoring with ${\mathcal O}(f)$ gives an isomorphism
$$\m_{v} \cong \m_{\tilde{v}}, \, \, \text{where} \, \, \tilde{v} = r + (\sigma + (k+r)f) + (p+1) \omega.$$ Such a twist raises the Euler characteristic by $1$. We normalize the moduli space by requiring that  $p = 1-r;$ when it has dimension $2a$ we refer to it as $\m_{r}^{a}$. Points in 
$\m_r^a$ are rank $r$ sheaves with Mukai vector 
$$v_{r,a} = r + (\sigma + (a - r(r-1))f)  + (1-r) \omega.$$ The normalization amounts to imposing that $$\chi (E) = 1 \, \, \text{for} \, \, E \in \m_{v}.$$ 
  
\vskip.1in

In rank $1$, note the isomorphism
$$X^{[a]} \cong \m_1^a, \, \, \, \, I_Z \mapsto I_Z (\sigma + af).$$ 
For any $r$, the generic point $E_r$ of $\m_r^a$ has exactly one section \cite{ogrady},
$$h^0 (E_r) = 1,$$
 as expected since the Euler characteristic is $1$.
Moreover, $h^0 (E_r (-f)) = 0$ generically, and 
$$h^0 (E_r (-2f)) = 0 \, \text{for} \, E_r \, \text{outside a codimension 2 locus in} \, \m_r^a.$$ In addition, stability forces the vanishing $h^2 (E_r (-2f)) = 0$ for all sheaves in $\m_r^a$, so
\begin{equation}
\label{h11}
h^1 (E_r (-2f)) = - \chi (E_r (-2f)) = 1
\end{equation}
 outside a codimension 2 locus in $\m_r^a$. In \cite{ogrady}, an open subscheme $U_r^a \subset \m_r^a$ is singled out, on which \eqref{h11} holds. For sheaves $E_r$ in $U_r^a$ there is a unique nontrivial extension
  \begin{equation}
\label{uniqueext}
0 \to {\mathcal O} \to {E}_{r+1}  \to E_r \otimes {\mathcal O} (-2f) \to 0.
\end{equation}
The resulting middle term ${E}_{r+1}$ is torsion-free, with Mukai vector $v_{r+1, a}$, and is stable unless $E_r$ belongs to a divisor 
$D_r$ in $U_r^a$. In the latter case, a stabilization procedure is required to ensure that the resulting rank $r+1$ sheaf also belongs to $\m_{r+1}^a$. 
The assignment $$E_r \mapsto E_{r+1}$$ identifies open subschemes $$U_r^a \cong U_{r+1}^a,$$ giving rise to a birational map
\begin{equation}
\label{isom}
\phi_r: \m_r^a \dashrightarrow \m_{r+1}^a, 
\end{equation}
and therefore a birational morphism away from codimension $2$,
 \begin{equation}
 \label{birational}
\Phi_r:  X^{[a]} \cong \m_1^a \dashrightarrow \m_r^a.
 \end{equation}

It will not be necessary for us to dwell on the details of the semistable reduction process along the $D_r$s although this, together with the identification of the $D_r$s themselves as divisors on the Hilbert scheme $X^{[a]}$, constitutes the most difficult part of \cite{ogrady}. We record here however, for future use, that 
\begin{equation}
\label{unstable}
D_1 = Q \cup S, \, \, \text{and}\, \, D_r = S\, \, \text{for} \, \, r \geq 2.
\end{equation}
Here, $Q$ is the divisor on $X^{[a]}$ consisting of ideals $I_Z$ such that $$h^0 (I_Z ((a-1)f)) \neq 0.$$ Equivalently, $Q$ is the divisor of cycles on $X^{[a]}$ with at least two points contained in the same elliptic fiber of $X$.  Furthermore, $S$ is the divisor of cycles in $X^{[a]}$ which intersect the section $\sigma$ of the elliptic fibration. 

\vskip.1in

\subsection {O'Grady's moduli space via Fourier-Mukai}\label{s22} We will reinterpret here the birational map $$\Phi_r:X^{[a]}\dasharrow \m_r^a$$ as a Fourier-Mukai transform. This will be crucial for the identification of the 
theta divisor and in particular for the proof of Theorem \ref{id}. 

We let $Y\to \mathbb P^1$ denote the dual elliptic $K3$ surface {\it i.e., } the relative moduli space of rank $1$ degree $0$ sheaves over the fibers of $\pi: X \to {\mathbb P}^1$. In fact, $X$ and $Y$ are canonically isomorphic. Writing $$\mathcal P\to X\times_{\mathbb P^1} Y$$ for the universal sheaf, we consider the Fourier-Mukai transform $$\mathsf S_{X\to Y}: {\mathbf D}(X)\to {\mathbf D}(Y), $$ with kernel $\mathcal P$, given by 
\begin{equation}
\mathsf S_{X\to Y} (x)=\mathbf {R}q_!\left( \mathcal P\stackrel{\mathbf L}{\otimes} {\mathbf L}p^{\star} x\right).
\label{fmfunctor}
\end{equation}
Here $p$ and $q$ are the two projections. We will normalize $\mathcal P$ such that $$\mathsf S_{X\to Y}(\mathcal O)=\mathcal O_{\sigma}[-1].$$ In fact, we have $$c_1(\mathcal P)=\Delta-p^{\star}\sigma-q^{\star}\sigma$$ where $\Delta$ is the diagonal in $X\times_{\mathbb P^1}Y$. In a similar fashion, we set $$\mathcal Q=\mathcal P^{\vee},$$ and use this as the kernel of the transform $$\mathsf T_{Y\to X}: \mathbf D(Y)\to \mathbf D(X)$$ It was shown in \cite {Br} that the functors $\mathsf S_{X\to Y}$ and $\mathsf T_{Y\to X}$ are equivalences of categories and that \begin {equation}\label{bridge}\mathsf S\circ \mathsf T=\mathbf 1_{\mathbf D(Y)}[-2],\,\,\,\, \mathsf T\circ \mathsf S=\mathbf 1_{\mathbf D(X)}[-2].\end{equation}

\vskip.2in

Fix a cycle $Z\in X^{[a]}$, and let $E_r$ denote the sheaf in $\m_r^a$ corresponding to $Z$ under the O'Grady isomorphism $\Phi_r$. We will consider {\it generic} subschemes $Z$, in the sense that 
\begin {itemize}
\item [(i)] $Z$ consists of distinct points,
\item [(ii)] no two points of $Z$ lie in the same fiber, 
\item [(iii)] $Z$ is disjoint from the section,
\item [(iv)] $Z$ does not contain any singular points of the fibers. 
\end {itemize} 

We determine the images of $E_r$ and of its derived dual $E_r^{\vee}$ under the functor \eqref{fmfunctor}. The answer is simpler for the dual, which is in fact $WIT_1$ relative to $\mathsf S_{X \to Y}$. Using the natural identification $Y\cong X$, we show 

\begin {proposition}\label{fmr} For generic $Z$, we have $$\mathsf S_{X\to Y}(E_r^{\vee})= I_Z(r\sigma+2rf)[-1].$$ Furthermore, $$\mathsf S_{X\to Y}(E_r^{\vee}(nf))= \mathsf S_{X\to Y}(E_r^{\vee})\otimes \mathcal O(nf).$$
\end {proposition}

The Fourier-Mukai transform of $E_r$ is expressed in terms of derived duals of ideal sheaves. Let $$\widetilde {Z}=\iota^{\star} Z$$ be the cycle obtained by taking the inverses of all points in $Z$ in the group law of their corresponding fibers. This makes sense even for singular fibers using the group law of the regular locus. We have

\begin {proposition} \label{fms}For $r\geq 1$, we have $$\mathsf S_{X\to Y}(E_r)=I_{\widetilde Z}^{\vee}\otimes \mathcal O(-r\sigma - 2(r-1) f).$$
\end {proposition}

\vskip.1in

The rest of this section is devoted to the proofs of Propositions \ref{fmr} and \ref{fms}. We study first how the generic sheaf $E_r$ restricts to the fibers. Consider a fiber $f$ of $\pi:X\to \mathbb P^1$ with origin $o=\sigma\cap f$, and let $$\mathsf W_r\to f$$ be the unique rank $r$ stable bundle on $f$ with determinant $\mathcal O_f(o)$. The $\mathsf W_r$'s were constructed by Atiyah over smooth elliptic curves. His arguments extend verbatim to nodal genus $1$ curves: we define $\mathsf W_r$ inductively as the unique nontrivial extension \begin{equation}\label{extension}0\to \mathcal O\to \mathsf W_{r+1}\to \mathsf W_{r} \to 0,\,\,\,\, \mathsf W_1=\mathcal O_f(o).\end{equation} 
Similarly, if $p$ is any smooth point of the fiber $f$, we write $$\mathsf W_{r, \,p}\to f$$ for the Atiyah bundle of determinant $\mathcal O_f(p)$ such that 
\begin{equation}
\label{extensionp}
0\to \mathcal O\to \mathsf W_{r+1,\, p}\to \mathsf W_{r, \,p}\to 0.
\end{equation}
The convention $$\mathsf W_{0, \,p}=\mathcal O_p$$ is used throughout. We have the following

\begin {lemma}\label{spl1} (i) If $f$ is a fiber such that $Z \cap f = \emptyset$, then $${E_r} {|_{f}} = \mathsf W_r.$$ 
(ii) If $f$ is a fiber through $p \in Z$, then $${E_r}{|_{f}} = \mathsf W_{r-1,\, p}\oplus \mathcal O_{f} (o-p).$$ 
\end {lemma}
\proof This is seen by induction starting with the case $r=1$ when $$E_1=I_Z(\sigma+af).$$ The basic observation is that for $p \in X$ and $I_p$ denoting its ideal sheaf in $X$, the restriction to the fiber $f$ through $p$ is $${I_p}{|_{f}} = {\mathcal O}_p \oplus {\mathcal O}_f (-p).$$ This gives the statement for $E_1$.
The inductive step from $r$ to $r+1$ follows from the exact sequence $$0\to \mathcal O\to E_{r+1} \to E_{r}(-2f)\to 0.$$ Its restriction to any fiber never splits as explained by Lemma $I.4.7$ \cite {ogrady}. The restriction to a fiber avoiding $Z$ must therefore coincide with the Atiyah bundle $\mathsf W_{r+1}$, since the latter is the only nontrivial extension $$0\to \mathcal O\to \mathsf W_{r+1}\to \mathsf W_{r}\to 0.$$ The same argument holds for the fibers through points of $Z$, using that there is a unique extension $$0\to \mathcal O\to \mathsf W_{r,\, p}\oplus \mathcal O_f(o-p)\to \mathsf W_{r-1, \,p}\oplus \mathcal O_f (o-p)\to 0.$$ \qed

Letting $f$ be a smooth elliptic fiber, we record now  the Fourier-Mukai transforms of the Atiyah bundles relative to the standard Poincar{\'{e}} kernel on $f \times f.$ We use hatted notation for the transforms, and as before we let $$\iota:f\to f$$ denote reflection about the origin of $f$. We have
\begin{equation}   
\label{zero}
\widehat {\mathsf W_r} = \mathcal O_f (-r\cdot o),
\end{equation}
\begin{equation}
\label{p}
\widehat {\mathsf W_{r, \,p}} = \mathcal O_f (-(r+1)\cdot o+\iota^{\star}p).
\end{equation}
By the results of \cite {muk}, the last two equations imply 
\begin{equation}
\label{zerodual}
\widehat{\mathsf W_r^{\vee}} = \mathcal O_f (r\cdot o)[-1],
\end{equation}
\begin{equation}
\label{pdual}
\widehat {\mathsf W_{r, \, p}^{\vee}} = \mathcal O_f((r+1)\cdot o-p)[-1].
\end{equation}

The first transform \eqref{zero} is obtained inductively by applying the Fourier-Mukai functor to the defining sequence \eqref{extension}. The base case $r=1$ is obvious. Similarly \eqref{p} can be derived using sequence \eqref{extensionp}. An alternate argument starts by noticing $$\mathsf W_{r, \,p} =\mathsf W_{r} \otimes M$$ where $$M^{r}=\mathcal O_f(p-o).$$ The line bundle $M$ corresponds to a point $m\in f$. Then $rm=p$ holds in the group law of the fiber. Using the properties of the Fourier-Mukai transform \cite {muk}, we obtain 
$$\widehat {\mathsf W_{r, \,p}}=t_{m}^{\star}\widehat {\mathsf W_{r}}=t_m^{\star}\mathcal O_f(-r\cdot o) =\mathcal O_f(-r\cdot m) =\mathcal O_f(-(r+1)\cdot o+\iota^{\star}p).$$

Equations \eqref{zero}, \eqref{p}, \eqref{zerodual}, and \eqref{pdual} also hold for the singular nodal fibers; this is explained by Lemma $2.13$, Definition $2.15$, and Remark $2.17$ in \cite{BK}. Note that the transforms in \cite{BK} are stated for the functor $\mathsf T_{Y\to X}$, but the results for the functor $\mathsf S_{X\to Y}$ follow via \eqref{bridge}. 

\vskip.2in

{\it Proof of Proposition \ref{fmr}.} We will first check that the isomorphism 
$$\mathsf S_{X\to Y}(E_r^{\vee})[1]=I_Z(r\sigma+2rf)$$ holds fiberwise. 
Derived restriction to fibers commutes with Fourier-Mukai \cite {Br}, and Lemma \ref{spl1} gives the restriction of $E_r^{\vee}$ to each fiber. The Fourier-Mukai transform of the restriction to a general fiber is 
$${\mathcal O}_f (r\cdot o)[-1],$$ by \eqref{zerodual}. For a special fiber $f$ containing a point $p \in Z,$ equation \eqref{pdual} yields the transform $$\mathcal O_f(r\cdot o-p)[-1]\oplus \mathcal O_p[-1].$$  The two formulas above give precisely the derived restriction of $I_Z(r\sigma+2rf)[-1]$. We have therefore checked that the proposition holds on every fiber. 

Since both sides are sheaves of rank $1$, we complete the proof by checking equality of determinants. As $X$ is simply connected, it is enough to match the first Chern classes. In general, let $V$ be a rank $r$ sheaf of Euler characteristic $\chi$ and $$c_1(V)=l\sigma+m f.$$ Then, by Grothendieck-Riemann-Roch, we have $$c_1(\mathsf S _{X\to Y}(V))=q_{!} \left (p^{\star} \text{ch}(V)\cdot \text {Todd}(X\times_{\mathbb P^1} Y/Y) \cdot \text {ch} \mathcal P\right)_{(2)}.$$ The Chern character of $V$ is $$\text{ch}(V)=r+(l\sigma+mf)+(\chi-2r)\omega,$$ where $\omega$ is the class of a point. Moreover, $$\text {Todd}(X\times_{\mathbb P^1} Y/Y)=p^{\star}(1-f+2\omega).$$ Hence $$c_1(\mathsf S_{X\to Y} (V))=rc_1(\mathsf S_{X\to Y}(\mathcal O))+(\chi -2r-l)q_{!}(p^{\star}\omega)+q_{!} \left(p^{\star} (l\sigma+mf)c_1(\mathcal P)\right)$$ $$=-r\sigma+(\chi-2r-l)f+2lf$$ $$=-r\sigma+(\chi-2r+l)f.$$ 

For $V=E_r^{\vee}$ the Chern class calculation gives $$c_1(\mathsf S_{X\to Y}(E_r^{\vee}))=-r\sigma-2rf,$$ which proves the first isomorphism. The calculation also shows that $$c_1(\mathsf S_{X\to Y}(E_r^{\vee}(nf))=c_1(\mathsf S_{X\to Y}(E_r^{\vee}))-nf.$$ The claim about twisting by fibers follows by repeating the above argument for $E_r^{\vee}(nf)$ and comparing determinants.\qed\vskip.1in

{\it Proof of Proposition \ref{fms}.} We consider the Fourier-Mukai functor $$\mathsf T_{X\to Y}:\mathbf D(X)\to \mathbf D(Y)$$ with kernel $$\mathcal Q=\mathcal P^{\vee}.$$ By duality, for all $x\in \mathbf D(X)$, we have $$\mathsf S_{X\to Y}(x)^{\vee}=\mathsf T_{X\to Y}(x^{\vee}\otimes \omega_{X/\mathbb P^1})[-1]=\mathsf T_{X\to Y}(x^{\vee}\otimes \mathcal O(2f))[-1].$$ Thus, the proposition follows once we establish that $$\mathsf T_{X\to Y}(E_r^{\vee}(2f))[-1]=I_{\widetilde Z}\otimes \mathcal O(r\sigma+2(r-1)f).$$ The proof of this fact is similar to that of Proposition \ref{fmr}. First, the equality is checked fiberwise using Lemma \ref{spl1}. The Grothendieck-Riemann-Roch calculation completes the argument. 

\vskip.1in

\begin {remark} The derived dual of the ideal sheaf $I_{\widetilde Z}$ can be computed explicitly for generic schemes $Z$. We include this calculation for completeness, even though it is not necessary for the proofs of the main theorems. 

First, O'Grady's construction gives rise to a rank $2$ sheaf $\widetilde{E_2}$  together with an exact sequence $$0\to \mathcal O\to\widetilde{E_2}\to I_{\widetilde Z}(\sigma+(a-2)f)\to 0.$$ Note that $\widetilde E_2$ is locally free by Lemma \ref{lf} below. Setting $$\mathsf C_Z=\left[\widetilde {E_2}\to \mathcal O(\sigma+(a-2)f)\right],$$ we claim \begin{equation}\label{dd}I_{\widetilde Z}^{\vee}=\mathsf C_{Z}.\end{equation} In particular, this implies that \begin{equation}\label{ddd}\mathsf S_{X\to Y}(E_r)=\mathsf C_Z \otimes \mathcal O(-r\sigma - 2(r-1) f).\end{equation} 

To prove \eqref{dd}, we dualize the sequence $$0\to \mathcal O\to \widetilde {E_2}\to I_{\widetilde Z}\otimes \mathcal O(\sigma+(a-2)f)\to 0.$$ We obtain $$0\to \mathcal O(-\sigma-(a-2)f)\to \widetilde{E_2}^{\vee}\to \mathcal O\to \mathcal {E}xt^1(I_{\widetilde Z}\otimes O(\sigma+(a-2)f), \mathcal O)\to 0.$$ It is well-known, see \cite {F} page $41$, that $$\mathcal {E}xt^1(I_{\widetilde Z}\otimes O(\sigma+(a-2)f), \mathcal O)=\mathcal O_{\widetilde Z}$$ hence the exact sequence yields $$0\to \mathcal O(-\sigma-(a-2)f)\to \widetilde{E_2}^{\vee}\to I_{\widetilde Z}\to 0.$$ Equation \eqref{dd} follows from here. \qed
\vskip.1in

\begin {lemma}\label{lf} If $Z$ contains no two points in the same fiber, then the associated sheaf $\widetilde E_2$ is locally free. 
\end {lemma} 

\proof Consider the divisor $Q$ of subschemes in $X^{[a]}$ containing two points in the same fiber. Let $$\mathcal D\hookrightarrow \m_2^a$$ be the codimension $1$ locus of nonlocally free sheaves in the rank $2$ moduli space. Lemma $4.41$ of \cite {Y} calculates $$\mathcal O(\mathcal D)=\Theta_w \, \, \, \text{on} \, \, \, \m_2^a,$$ for the Mukai  vector $$w=(2, -\sigma - (a-2)f, (a-2)\omega).$$ Using now the formulas in \cite {ogrady}, page $27$, and \eqref{qaa} below, we have
 $$\Theta_w=\mathcal O(Q)$$ under the identification $$X^{[a]}\dasharrow \m_2^a.$$ Finally, we will remark in \eqref{ss} below that the line bundle $\mathcal O(Q)\to X^{[a]}$ has a unique section, hence $\mathcal D$ and $Q$ coincide as claimed. \qed
 
\end {remark}

\subsection {Line bundles and theta divisors over the Hilbert scheme of points} The birational isomorphism \eqref{birational} allows us to identify the Picard group of $\m_r^a$ with that of the Hilbert scheme $X^{[a]}$. 

For any smooth projective surface $X$ and any line bundle $L$ on it, we indicate by $L_{(a)}$ the line bundle on $X^{[a]}$ induced from the symmetric line bundle $L^{\boxtimes a}$ on the product $X\times \ldots \times X$.  
Letting $p$ and $q$ be the projections $$p:X^{[a]} \times X\to X^{[a]}, \,\,\,q:X^{[a]}\times X\to X,$$ and letting ${\mathcal O}_{\mathcal Z}$ denote the universal structure sheaf on 
$X^{[a]} \times X$, we further set 
\begin{equation}
\label{pic}
L^{[a]} = \det \, p_{!} ({\mathcal O}_{\mathcal Z} \otimes q^{\star} L).
\end{equation} 
It is well known that the line bundles $L_{(a)}$ for $L \in {\text{Pic}}\, X$ and $M = {\mathcal O}^{[a]}$ generate the Picard group of $X^{[a]}$, and that for any $L$ on $X$, $$L^{[a]} = L_{(a)} \otimes M.$$ We have, for instance, 
$${\mathcal O}(S) = {\mathcal O} (\sigma)_{(a)},$$ and \begin{equation}\label{qaa}\mathcal O(Q)={\mathcal O ((a-1)f)}^{[a]}.\end{equation} We note for future use the formulas of \cite {egl}, \begin{equation}\label{sections}h^0(X^{[a]}, L_{(a)})=\binom{h^{0}(X, L)+a-1}{a},\,\,\,\,  h^0(X	^{[a]}, L^{[a]})=\binom{h^0(X,L)}{a}.\end {equation} To illustrate, we compute \begin {equation}\label{ss}h^0(X^{[a]}, \mathcal O(Q))= h^0(X^{[a]}, \mathcal O((a-1)f)^{[a]})=\binom {h^0(\mathcal O((a-1)f))}{a}=1.\end {equation} 

\vskip.1in

Consider now two Hilbert schemes of points $X^{[a]}$ and $X^{[b]}$, and the rational morphism, defined away from codimension 2, 
\begin {equation}
\tau: X^{[a]} \times X^{[b]} \dashrightarrow X^{[a+b]}, \, \, \, \, \, (I_Z, I_W) \mapsto I_Z \otimes I_W.
\end{equation} Assume that $L$ is a line bundle on $X$ with no higher cohomology, and such that 
$$\chi(L) = h^0 (L) = a+b.$$
From \eqref{sections}, we note that
 $$h^0 (X^{[a+b]}, L^{[a+b]}) = \binom{h^0 (X, L)}{a+b} = 1.$$ The unique section of $L^{[a+b]}$ vanishes on the locus
 \begin{equation}
\label{b-theta}
 \theta_L = \{I_V \in X^{[a+b]}, \, \, \,\text{such that} \, \, \, H^0 (I_V \otimes L) \neq 0\},
 \end{equation}
 whose pullback under $\tau$ is the 
 divisor $$\theta_{L,a,b} = \{(I_Z, I_W) \in X^{[a]} \times X^{[b]} \, \, \, \text{such that} \, \, \, H^0 (I_Z \otimes I_W \otimes L) \neq 0\}.$$ 
We furthermore have
\begin{equation}
\mathcal O (\theta_{L, a, b}) = \tau^{\star} L^{[a+b]} = L^{[a]} \boxtimes L^{[b]} \, \, \, \text{on} \, \, \, X^{[a]} \times X^{[b]}. 
\end{equation}
It was observed in \cite{mo} that $\theta_{L,a,b}$ induces an isomorphism 
\begin{equation}
{\mathsf D}: H^0 (X^{[a]}, L^{[a]})^{\vee} \to H^0 (X^{[b]}, L^{[b]}).
\end{equation} 

It will be important for our arguments to consider the following partial reflection of the divisor $\theta_{L, a, b}$: 
$$\widetilde \theta_{L, a, b}=\{(Z, W) \in X^{[a]} \times X^{[b]} \text { such that } h^0(I_Z\otimes I_{\widetilde W}\otimes L)\neq 0\}.$$ As usual, the subschemes $$\widetilde Z=\iota^{\star} Z, \,\, \widetilde W=\iota^{\star} W$$ are obtained from the fiberwise reflection $\iota:X\dasharrow X$. There is a seeming asymmetry in the roles of $Z$ and $W$ in the definition of $\widetilde \theta_{L, a, b}$, but in fact we also have $$\widetilde \theta_{L, a, b}= \{(Z, W) \in X^{[a]} \times X^{[b]} \text { such that } h^0(I_{\widetilde Z}\otimes I_W\otimes L)\neq 0\}.$$ To explain this equality, note first that the line bundle $L$ is invariant under $\iota$ $$\iota^{\star}L=L.$$ Hence, so are the tautological line bundles $L^{[a]}$, $L^{[b]}$ and $L^{[a+b]}$. On $X^{[a+b]},$ the divisor $\theta_L$ of \eqref{b-theta} corresponds to the unique section of $L^{[a+b]},$ therefore must be invariant under $\iota$ as well, 
$$\iota^{\star} \theta_L = \theta_L.$$ The same is then true for the pullback $$\theta_{L,a,b} = \tau^{\star} \theta_L,$$ which implies that
$$h^0 (I_Z \otimes I_{\widetilde{W}} \otimes L) = 0 \, \, \, \text{if and only if} \, \, \, h^0 (I_{\widetilde{Z}} \otimes I_{{W}} \otimes L) = 0.$$

The above discussion also shows that $\tilde \theta_{L, a, b}$ is a section of the line bundle $L^{[a]}\boxtimes L^{[b]}$ and that furthermore it induces an isomophism \begin{equation}\label{dt} \widetilde {\mathsf D}: H^0 (X^{[a]}, L^{[a]})^{\vee} \to H^0 (X^{[b]}, L^{[b]}).
\end{equation} 
 
 \vskip.1in
 
\subsection{The strange duality setup and the standard theta divisor}
We now place ourselves in the setting of Theorems \ref{sd}, \ref{t3} and \ref{id} {\it{i.e.,}} we take $X$ to be an elliptically fibered $K3$ surface with section, and consider two moduli spaces of sheaves $\m_v$ and $\m_w$ with orthogonal Mukai vectors satisfying \begin{itemize}\item [(i)] $\langle v, w^{\vee} \rangle=0$, \item [(ii)]$c_{1}(v)\cdot f=c_{1}(w)\cdot f=1$, \item [(iii)] $\langle v, v\rangle + \langle w, w\rangle \geq 2(r+s)^{2}.$\end{itemize} 
Equivalently, we consider two normalized moduli spaces $\m_r^a$ and $\m_s^b$ such that
\begin{equation}
\label{divisibility}
r+s \, | \, a + b - 2, \, \, \text{and moreover} \, \, -\nu \, =_{\text{def}}\, \frac{a+b -2}{r+s} - (r+s -2)  \geq 2. 
\end{equation}
We also assume that $r, s \geq 2.$ The divisibility condition and the definition of $\nu$ are so as to ensure that 
$$\chi (E_r \cdot F_s  (\nu f)) = 0, \, \, \text{for sheaves} \, \, E_r \in \m_r^a, F_s \in \m_s^b.$$
Furthermore, the stability condition implies that 
$$H^2 (E_r \otimes F_s (\nu f)) = 0.$$ The vanishing $$\text {Tor}^1(E_r, F_s)=\text{Tor}^2(E_r, F_s)=0$$ is satisfied when $E_r$ or $F_s$ are locally free, which occurs away from codimension $2$ in the product space.

We denote by $\Theta_{r,s}$ the locus 
$$\{ (E_r, F_s) \in \m_r^a \times \m_s^b \, \, \text{such that} \, \, h^0 (E_r \otimes F_s (\nu f)) \neq 0 \}.$$ The condition defining $\Theta_{r,s}$ is divisorial, but it is not a priori clear that this locus actually has codimension $1$. Nonetheless, using the explicit formulas of \cite{ogrady}, the {\em{line bundle}} ${\mathcal O} (\Theta_{r,s})$ on $\m_r^a \times \m_s^b$ can be expressed on the product $X^{[a]} \times X^{[b]}$ via the birational map 
$$(\Phi_r, \Phi_s): X^{[a]} \times X^{[b]} \dashrightarrow \m_r^a  \times \m_s^b.$$
Letting
\begin{equation}
L = {\mathcal O} ((r+s)\sigma + (2 (r+s) - 2- \nu)f) \, \, \, \text{on} \, \, \, X,
\end{equation}
it was shown in \cite{mo} that 
\begin{equation}
{\mathcal O} (\Theta_{r,s}) = L^{[a]} \boxtimes  L^{[b]}.
\end{equation}
We prove that
\addtocounter{theorem}{-2}
\begin{theorem}
\label{id}
${\Theta_{r,s}} = \tilde{\theta}_{L,a,b}$ on $X^{[a]} \times X^{[b]}.$
\end{theorem}
\vskip.1in

\subsection {The theta divisor over the generic locus} In this section and the one following it, we prove Theorems \ref{sd}, \ref{t3} and \ref{id}. 

We first identify the theta divisor $\Theta_{r,s}$ on the locus corresponding to generic $Z$ and $W$. Our genericity assumptions were specified in (i)-(iv) of Section \ref{s22}. On any Hilbert scheme of points of $X$, we consider then the following:
\begin {itemize}
\item [(i)] the divisor $R$ consisting of cycles with at least two coincident points,
\item [(ii)] the divisor $Q$ of cycles with two points on the same fiber,
\item [(iii)] the divisor $S$ of cycles which intersect the section.
\end {itemize} Recall that along the divisors $S$ and $Q$ the extensions \eqref{uniqueext} have unstable middle terms needing to undergo semi-stable reduction.

We single out here only the nongeneric loci corresponding to divisors, as for our purposes we can ignore higher codimension phenomena. Thus we will disregard the loci corresponding to 
\begin {itemize}
\item [(iv)] schemes whose supports contain singular points of fibers of $X.$ 
\end {itemize}

We work with the rational morphism $$\tau: X^{[a]} \times X^{[b]} \dashrightarrow  X^{[a+b]},$$ and we will pullback the divisors $R$, $S$ and $Q$ to the product of Hilbert schemes and of moduli spaces $\m_r^a$ and $\m_s^b$. We set $$\mathfrak M=\m_r^a\times \m_s^b\setminus (\tau^{\star}R\cup \tau^{\star}Q\cup \tau^{\star} S).$$ 
For $(E_r, F_s)\in \m$, we show
\begin{equation}
\label{parseval}
H^0 (E_r \otimes F_s (\nu f)) = 0 \, \, \, \text{if and only if} \, \, \, H^1 (I_Z \otimes I_{\widetilde{W}} \otimes L) = 0.
\end{equation}
In other words, we prove \begin{equation}
\label{generic}
\Theta_{r,s}  \setminus (\tau^{\star}R\cup \tau^{\star}Q\cup \tau^{\star} S) = \tilde{\theta}_{L,a,b} \setminus (\tau^{\star}R\cup \tau^{\star}Q\cup \tau^{\star} S).
\end{equation}

\vskip.1in

To establish \eqref{parseval}, we use the Fourier-Mukai functor ${\mathsf S}_{X \to Y}$ defined in \eqref{fmfunctor}, as well as Propositions 1 and 2. 
We calculate \begin {eqnarray*} H^0(E_r\otimes F_s(\nu f))&=&\text {Hom}_{\mathbf D(X)} (E_r^{\vee}(-\nu f), F_s)\\ &=&\text{ Hom}_{\mathbf D(Y)} \left(\mathsf S_{X\to Y} (E_r^{\vee}(-\nu f)), \mathsf S_{X\to Y} (F_s)\right)\\&=&\text{ Hom}_{\mathbf D(Y)} \left( I_Z(r\sigma+2rf-\nu f), I_{\widetilde W}^{\vee} \otimes \mathcal O(-s\sigma -2(s-1) f)[1]\right)\\&=&\text{Ext}^{1}(I_Z\otimes L,\, I_{\widetilde W}^{\vee})\\&=&\text {Ext}^{1} (I_{\widetilde W}^{\vee},\, I_Z\otimes L)^{\vee}\\&=&H^{1}(I_{\widetilde W}\otimes I_Z\otimes L)^{\vee}.\end{eqnarray*} This calculation renders \eqref{parseval} obvious.

Note that Theorem \ref{t3} is a consequence of the above discussion. \qed
\vskip.1in

\subsection {The nongeneric locus} 

In order to complete the proof of Theorem \ref{id}, we will need to analyze the overlaps of the theta divisor with $R$, $Q$ and $S$. 

First, as the divisors $R, Q, S$ are invariant under $\iota,$ we write \eqref{generic} equivalently as 
\begin{equation}
\label{generic2}
\widetilde{\Theta}_{r,s}  \setminus (\tau^{\star}R\cup \tau^{\star}Q\cup \tau^{\star} S) = {\theta}_{L,a,b} \setminus (\tau^{\star}R\cup \tau^{\star}Q\cup \tau^{\star} S),
\end{equation} where $\widetilde \Theta_{r, s}$ is the partial reflection of the divisor $\Theta_{r, s}$ obtained by acting with the involution $\iota$ on one of the factors. 

We write
$$\widetilde{\Theta}_{r,s} = \Gamma \cup \Delta,$$ where $\Gamma$ and $\Delta$ are divisors such that the intersection 
$$\Delta\cap  \left ( \tau^{\star} Q\cup \tau^{\star}R \cup \tau^{\star} S\right )$$ is proper, and $$\text{ support } \Gamma \subset \tau^{\star} Q \cup \tau^{\star} R\cup \tau^{\star} S.$$
Equation \eqref{generic2} shows in particular that $\Delta$ is a pullback divisor under $\tau,$ $$\Delta = \tau^{\star} {\Delta_0} \, \, \text{for} \, \, \Delta_0 \subset X^{[a+b]}.$$ Since $${\mathcal O} (\widetilde{\Theta}_{r,s}) = L^{[a]} \boxtimes L^{[b]} = \tau^{\star} L^{[a+b]},$$ we have \begin{equation}\label{gammaeq}{\mathcal O} (\Gamma) = \tau^{\star} (L^{[a+b]} \otimes {\mathcal O} (-\Delta_0)).\end{equation} More strongly, we will show shortly that \eqref{gammaeq} implies that 
\begin{claim} \label{cl1}$\Gamma$ as a divisor is a pullback under the morphism $\tau,$ $$\Gamma=\tau^{\star} \Gamma_0.$$ \end{claim}

As a consequence, $$ \widetilde{\Theta}_{r,s} = \tau^{\star} (\Delta_0 \cup \Gamma_0).$$ Now $$\theta_L = \{V \, \, \text{such that} \, \, h^0 (I_V \otimes L) \neq 0 \} $$ is the only section of $L^{[a+b]}$ on $X^{[a+b]}$. Thus we must have that
$$\theta_L = \Delta_0 \cup \Gamma_0, \, \, \, \text{and}$$ $$\widetilde{\Theta}_{r,s} = \tau^{\star} \theta_L = \{ ( I_Z, I_W) \, \, \text{such that} \, \, h^0 (I_Z \otimes I_W \otimes L ) \neq 0\} = \theta_{L,a,b}.$$ This completes the proof of Theorem \ref{id}. Theorem \ref{sd} follows as well via \eqref{dt}. \qed

\vskip.1in

{\it Proof of Claim \ref{cl1}.} We will consider the three divisors $Q$, $R$ and $S$ over the Hilbert schemes $X^{[a]},$ $X^{[b]}$ or $X^{[a+b]}.$ All these divisors are irreducible. Let $$\tau^{\star}Q=Q_1\cup Q_2\cup Q_3,\,\,\,\,\,\,\tau^{\star}R=R_1\cup R_2,$$ $$\tau^{\star}S=S_1\cup S_2$$ be the irreducible components of the pullbacks on the product $X^{[a]}\times X^{[b]}$. Here $$Q_1=Q\times X^{[b]},\,\,\,\,Q_2=X^{[a]}\times Q,$$ while $Q_3$ is the divisor of cycles $(I_Z, I_W)\in X^{[a]}\times X^{[b]}$ such that $Z, W$ intersect the same elliptic fiber. In the same fashion, we have $$R_1=R\times X^{[b]}, R_2=X^{[a]}\times R,$$ $$S_1=S\times X^{[b]}, S_2=X^{[a]}\times S.$$ 

Note first that $\widetilde{\Theta}_{r,s}$ does not contain the divisor $Q_3$. Indeed, $$\mathcal O(Q)=\mathcal O((a+b-1)f)^{[a+b]} \text { on } X^{[a+b]},$$ so $$\tau^{\star}\mathcal O(Q)=\mathcal O((a+b-1)f)^{[a]}\boxtimes \mathcal O((a+b-1)f)^{[b]} \text { on } X^{[a]}\times X^{[b]}.$$ As $$\mathcal O(Q_1)=\mathcal O((a-1)f)^{[a]}\boxtimes \mathcal O\,\,\, \text { and }\,\,\, \mathcal O(Q_2)=\mathcal O\boxtimes \mathcal O((b-1)f)^{[b]},$$ it follows that $$\mathcal O(Q_3)=\mathcal O(bf)_{(a)}\boxtimes \mathcal O(af)_{(b)}.$$ Assuming $\widetilde{\Theta}_{r,s}$ contained $Q_3$, we would have \begin{equation}\label{cohgr}H^{0}(X^{[a]}\times X^{[b]}, \mathcal O(\widetilde{\Theta}_{r,s}-Q_3))\neq 0.\end{equation} However, we will show that \eqref{cohgr} is false. Indeed, $$\mathcal O(\widetilde{\Theta}_{r,s}-Q_3)=L(-bf)^{[a]}\boxtimes L(-af)^{[b]}.$$ From \eqref{sections}, we have $$h^{0}\left(L(-bf)^{[a]}\right)=\binom{h^0(L(-bf))}{a}, $$ $$h^0\left(L(-af)^{[b]}\right)=\binom{h^0(L(-af))}{b}.$$ It suffices to explain that either \begin{equation}\label{h0}h^{0}(L(-bf))=0 \text { or } h^0(L(-af))=0.\end{equation}

On the surface $X$, we generally have \begin{equation} \label{surface} h^0 (X, {\mathcal O} (m\sigma + nf )) = \left \{ \begin{array} {l} 0, \, \, \, \text{if} \, \,\,  m \geq 0, \, \, n <0 \\ 2 + m (n-m), \, \, \, \text{if} \, \, \, m > 0, \, \, n \geq 2m \end{array} \right .  .\end{equation} The first dimension count is immediate as $$h^0 (X, {\mathcal O} (m\sigma)) = 1$$ for all $m \geq 0$, and the second holds as in that case ${\mathcal O} (m\sigma + n f)$ is big and nef, so has no higher cohomology. Now, recall that $$L=\mathcal O\left((r+s)\sigma+\left(r+s+\frac{a+b-2}{r+s}\right)f\right) \text { on } X.$$ The numerical constraint \eqref{divisibility} $$a+b\geq (r+s)^2+2$$ ensures that either $L(-af)$ or $L(-bf)$ has a negative number of fiber classes. This proves \eqref{h0} using the dimension count \eqref{surface}.

Similarly, $\widetilde{\Theta}_{r,s}$ cannot contain both $Q_1$ and $Q_2$. Indeed, we calculate $$\mathcal O(\widetilde{\Theta}_{r,s}-Q_1-Q_2)= L((-a+1)f)_{(a)}\boxtimes L((-b+1)f)_{(b)}.$$ As in \eqref{h0}, unless $$r=s=2, \, a=b=9,$$ we have \begin{equation}\label{h00}h^0(L((-a+1)f))=0 \text { or } h^0(L((-b+1)f))=0,\end{equation} therefore \eqref{sections} implies that $\mathcal O(\widetilde{\Theta}_{r,s}-Q_1-Q_2)$ has no sections. 

Let us write $$\Gamma=q_1 Q_1+q_2 Q_2+q_3Q_3+r_1 R_1+r_2R_2+s_1 S_1+s_2 S_2.$$ The above argument shows that $q_3=0$ and that we can assume without loss of generality $q_2=0$. We calculate $$\mathcal O(Q_1)=\mathcal O((a-1)f)^{[a]}\boxtimes \mathcal O,$$ $$\mathcal O(R_1)=M^{-2}\boxtimes \mathcal O, \,\,\,\mathcal O(R_2)=\mathcal O\boxtimes M^{-2},$$ $$\mathcal O(S_1)=\mathcal O(\sigma)_{(a)}\boxtimes \mathcal O, \,\,\, \mathcal O(S_2)=\mathcal O\boxtimes \mathcal O(\sigma)_{(b)},$$ 
Consequently, 
\begin{equation}\label{odegamma}\mathcal O(\Gamma)=\left(\mathcal O\left(q_1(a-1)f+s_1\sigma\right)_{(a)}\otimes M^{q_1-2r_1}\right)\boxtimes \left(\mathcal O\left(s_2\sigma\right)_{(b)}\otimes M^{-2r_2}\right).\end{equation} From \eqref{gammaeq} we know that this line bundle is a pullback under $\tau$. This strongly constrains the coefficients in the expression \eqref{odegamma}. In fact, via the isomorphism $$\text {Pic}(X^{[n]})=\text {Pic}(X)\oplus \mathbb Z,\,\,\,\,\, L_{(n)}\otimes M^{r}\mapsto (L, r),$$ the image $$\tau^{\star}:\text {Pic}(X^{[a+b]})\to \text {Pic}(X^{[a]})\times \text {Pic}(X^{[b]})$$ corresponds to the diagonal embedding. Therefore in \eqref{odegamma} we must have $$q_1=0, s_1=s_2, q_1-2r_1=-2r_2.$$ Hence $$\Gamma=r_1(R_1+R_2)+s_1(S_1+S_2)$$ is a pullback of the divisor $\Gamma_0=r_1R+s_1S.$  This establishes Claim \ref{cl1} unless $$r=s=2, \, \, a=b=9.$$

In this latter case, the dimension calculated in \eqref{h00} is $$h^0(L((-a+1)f))=h^0(L((-b+1)f))=h^0(\mathcal O(4\sigma))=1,$$  and $$\mathcal O(\widetilde \Theta_{r,s}-Q_1-Q_2)={\mathcal O} (4\sigma)_{(a)}\boxtimes {\mathcal O} (4\sigma )_{(b)}$$ has a unique section supported on $S_1\cup S_2$. Thus if $Q_1, Q_2$ are both contained in $\widetilde \Theta_{r,s}$, then $\widetilde \Theta_{r,s}$ is supported in $Q_1\cup Q_2\cup S_1\cup S_2$. This implies via \eqref{generic2} that $$\theta_L\setminus (R\cup Q\cup S)=\emptyset.$$ However, the following remark shows that this is not possible. 
\qed

\begin{remark}\label{r2}
We note here that in fact $\theta_L$ on $X^{[a+b]}$ intersects both $S$ and $Q$ properly, therefore $\widetilde{\Theta}_{r,s}$ intersects $\tau^{\star}S$ and $\tau^{\star} Q$ properly.
Otherwise, we would have
$$H^0 (X^{[a+b]}, L^{[a+b]} \otimes {\mathcal O} (-S)) \neq 0, $$ $$H^0 (X^{[a+b]}, L^{[a+b]} \otimes {\mathcal O} (-Q)) \neq 0.$$
We calculate $$L^{[a+b]}\otimes \mathcal O(-S)=L(-\sigma)^{[a+b]}.$$ From the dimension count \eqref{surface} we have $$h^0(L(-\sigma))=a+b+1+\nu,$$ and therefore from \eqref{sections}, $$h^{0}(L^{[a+b]}\otimes \mathcal O(-S))=\binom{h^{0}(L(-\sigma))}{a+b}=\binom{a+b +1 +\nu }{a+b} = 0.$$ 
Similarly, $$L^{[a+b]}\otimes \mathcal O(-Q) = L ((-a-b+1) f)_{(a+b)}$$ has no sections on $X^{[a+b]}$ since $L ((-a-b+1) f)$ has no sections on $X$.

\end {remark}

\vskip.1in

\subsection {Arbitrary elliptic surfaces}

Theorem \ref{ellsur} will be proved in this section. First, we write down O'Grady's construction for arbitrary simply connected elliptic surfaces with section, and then we reinterpret it via Fourier-Mukai transforms. 

The holomorphic Euler characteristic of the fibration $$\pi:X\to \mathbb P^1$$ will be denoted $$\chi=\chi(\mathcal O)=1+h^2(\mathcal O_X)>0.$$ 
We study normalized moduli spaces of sheaves $\m_v$ such that $$\chi(v)= 1 \implies c_1(v)=\sigma+\left(a-\frac{r(r-1)}{2}\chi\right)f,$$ where we write $2a$ for the dimension of $\m_v.$ A birational isomorphism $$\Phi_v:X^{[a]}\dasharrow \m_v$$ is constructed as follows. As in the case of $K3$ surfaces, we consider generic schemes $Z$ of length $a$, satisfying the requirements (i)-(iv) of section \ref{s22}. 
We set $$E_1=I_Z(\sigma+af).$$ Inductively, we construct nontrivial extensions \begin{equation}\label{ogrgen} 0\to \mathcal O\to E_{r+1}\to E_r(-\chi f)\to 0.\end{equation} Several statements are to be proved simultaneously during the induction step:
\begin {itemize}
\item [(a)] $\text{Ext}^0(E_r(-\chi f), \mathcal O)=0$
\item [(b)] $\text{Ext}^2(E_r(-\chi f), \mathcal O)=0.$
\item [(c)] $\text{Ext}^1(E_r(-\chi f), \mathcal O)\cong \mathbb C.$  This shows that the extension \eqref{ogrgen} is unique. 
\item [(d)] the restriction of $E_r$ to the generic fiber is the Atiyah bundle $\mathsf W_r$. This implies the stability of $E_r$ with respect to suitable polarizations. For special fibers through $p\in Z$, the restriction splits as $\mathsf W_{r-1,p}\oplus \mathcal O_f(o-p).$ 
\end {itemize} 
Checking (a)-(d) for the base case $r=1$ uses the requirements (i)-(iv) of section $2.2$. We briefly explain the inductive step. The first vanishing in (a) follows by stability since for polarizations $H=\sigma+mf$ with $m>>0$, we have $$\frac{c_1(E_r(-\chi f))\cdot H}{r}>0.$$ Regarding (b), we consider the exact sequence induced by \eqref{ogrgen} $$\text{Ext}^2(E_r(-2\chi f), \mathcal O)\to \text{Ext}^2(E_{r+1}(-\chi f), \mathcal O)\to \text{Ext}^2(\mathcal O(-\chi f), \mathcal O)=0.$$ Now (b) follows since the leftmost term also vanishes, as one can see by considering the injection $$\text{Ext}^2(E_{r}(-2\chi f), \mathcal O)\hookrightarrow \text{Ext}^2(E_{r}(-\chi f), \mathcal O)=0.$$ Now, (a) and (b) imply (c) via a Riemann-Roch calculation. Finally, (d) is argued exactly as Lemma \ref{spl1} above. 

We use (d) to calculate Fourier-Mukai transforms $$\mathsf S_{X\to Y}:\mathbf D(X)\to \mathbf D(Y).$$ By the arguments of section \ref{s22}, we obtain
\begin {itemize}
\item [(i)] $\mathsf S_{X\to Y}(E_r^{\vee})=I_Z(r\sigma+r\chi f)[-1]$
\item [(ii)] $\mathsf S_{X\to Y}(E_r)=I_{\widetilde Z}^{\vee} \otimes \mathcal O(-r\sigma-(r-1)\chi f).$
\end {itemize} 

Consider now two complementary moduli spaces $\m_v$ and $\m_w$. After twisting by fiber classes, we may assume $v$ and $w$ are normalized. Consider the theta locus $$\Theta=\{(E, F)\in \m_v\times \m_w:h^0(E\otimes F\otimes \mathcal O(\nu f))=0\}$$ where $$-\nu =\frac{a+b-\chi}{r+s} -(r+s-1)\frac{\chi}{2}+1\geq \chi,$$ by the condition (ii) of Theorem \ref{ellsur}. We set $$L=\mathcal O_X\left((r+s)\sigma+((r+s-1)\chi-\nu )f\right)\otimes K_X.$$ An easy calculation shows $$h^0(L)=\chi (L)=a+b,$$ and that $L$ has no higher cohomology. We therefore obtain a divisor $$\theta_{L, a, b}\subset X^{[a]}\times X^{[b]},$$ and the associated twist $$\tilde \theta_L=(1\times i)^{\star} \theta_{L}=(i\times 1)^{\star} \theta_L$$ in the product of Hilbert schemes. 

Repeating the argument for elliptic $K3$s, we obtain that under the birational map $$\Phi_v\times \Phi_w:X^{[a]}\times X^{[b]}\dasharrow \m_v\times \m_w$$ we have \begin{equation}\label{pvw}(\Phi_v\times \Phi_w)^{\star} \Theta \subset \tilde \theta_{L},\end{equation} at least along the nongeneric locus. This is enough to establish that $\Theta$ is a divisor, proving Theorem \ref{ellsur}. \qed

\begin {remark} Unfortunately, we cannot conclude equality in \eqref{pvw} since we are unable to estimate the codimension of the image of $\Phi_v$ and $\Phi_w$ in the two moduli spaces. We believe this codimension to be at least $2$. Conjecture \ref{ellsurconj} follows once this statement is established. 
\end{remark}

\section {Generic strange duality} In this section we prove Theorem \ref{genericarb} by a deformation argument. We will keep the same notations as in the introduction. 

Let $(X, H)$ be a polarized $K3$ surface, and let $$v=(r, H, \chi-r), \,\,\,w=(s, H, \chi'-s)$$ be the two orthogonal Mukai vectors with $\chi \leq 0, \chi'\leq 0$. Consider a deformation $$\pi:(\mathcal X, \mathcal L)\to \Delta$$ of polarized $K3$ surfaces such that \begin {itemize}
\item [(i)] the $K3$ surface $(X, H)$ appears as a generic fiber. We may assume that for $t\neq 0$, $\mathcal L_t$ is the unique ample generator of the Picard group of $\mathcal X_t$.  
\item [(ii)] $\mathcal X_0$ is an elliptically fibered $K3$ surface with a section, such that $c_1(\mathcal L_0)$ is a numerical section.  
\end {itemize}

For each $t\in \Delta$, we consider the Mukai vectors $$v_t=(r, c_1(\mathcal L_t), \chi-r), \,\,\, w_t=(s, c_1(\mathcal L_t), \chi'-s).$$ Since intersection products are preserved by deformations, we have $$\langle v_t^{\vee}, w_t\rangle =0 \text { for all } t\in \Delta.$$ We form two relative moduli spaces of $\mathcal L_t$-semistable sheaves $$\mathfrak M [v]=\cup_{t\in \Delta}\mathfrak M_{v_t}\to \Delta,\,\,\,\, \mathfrak M[w]=\cup_{t\in \Delta}\mathfrak M_{w_t}\to \Delta.$$ The product $$\mathfrak M[v]\times_{\Delta} \mathfrak M[w]\to \Delta$$ carries a relative theta divisor $\Theta[v,w]$ obtained as the vanishing locus of a section of the relative theta bundles $$\Theta[w]\boxtimes \Theta[v]\to \m[v]\times_{\Delta} \m[w].$$ Pushing forward to $\Delta$ via the projections $\pi$, we obtain the sheaves $$\mathsf V=\pi_{\star}\left(\Theta[w]\right),\,\, \mathsf W=\pi_{\star}\left(\Theta[v]\right),$$ as well as a section $\mathsf D$ of $\mathsf V\otimes \mathsf W.$ 

We claim that the sheaves $\mathsf V$ and $\mathsf W$ are locally free of equal rank. Let us consider $\mathsf V$ first. Over the special fiber, there is a birational isomorphism $$\m_{v_0}\dasharrow \mathcal X_0^{[a]}$$ regular away from codimension $2$, where $a$ denotes half the dimension of $\m_{v_0}$. Similarly $2b$ is the dimension of $\m_{w_0}.$ The line bundle $\Theta_{w_0}$ corresponds to a line bundle of the form $L^{[a]}$ for some $L\to \mathcal X_0$ with $h^0(L)=a+b.$ Hence, by \eqref{sections}, $$h^0(\m_{v_0}, \Theta_{w_0})=h^0(\mathcal X_0^{[a]}, L^{[a]})=\binom {h^0(L)}{a}=\binom {a+b}{a}.$$ Over the general fiber, Lemma \ref{nohc} below shows that $$h^0(\m_{v_t}, \Theta_{w_t})=\chi(\m_{v_t}, \Theta_{w_t})=\binom{a+b}{a}.$$ The calculation of the Euler characteristic in the equation above can be found in \cite {OG2}. By Grauert's theorem, $\mathsf V$ is a vector bundle whose formation commutes with restriction to fibers. The same arguments apply to $\mathsf W$. 

\begin {lemma} \label{nohc} Let $H$ be an ample generator of the Picard group of the $K3$ surface $X$. Assume that $v$ and $w$ are Mukai vectors such that \begin {itemize}
\item [(i)]$\langle v^{\vee}, w\rangle =0,$
\item [(ii)] $c_1(v)=c_1(w)=H,$
\item [(iii)] $\chi(v)\leq 0,\,\, \chi(w)\leq 0.$
\end {itemize}
The line bundle $\Theta_w\to \m_v$ is big and nef, hence it does not have higher cohomology. 
\end {lemma}

\proof For a Mukai vector $v=(v_0, v_2, v_4)$, define $$\lambda_v=(0, -v_0 H , H\cdot v_2)$$ and $$\mu_v=(-H\cdot v_2, v_4H, 0).$$ These vectors are perpendicular to $v$. It was shown by Jun Li that $\Theta_{-\lambda_v}$ is big and nef \cite {jli}; in fact, $\Theta_{-\lambda_v}$ defines a morphism from the Gieseker to the Uhlenbeck compactification. 

Using reflections along rigid sheaves, Yoshioka proved that $\Theta_{-\lambda_v-\mu_v}$ is also big and nef \cite {Y}, and that it determines a morphism $$\pi: \m_v\to \mathfrak X,$$ where $$\mathfrak X\subset \bigcup_{k\geq -\chi(v)} \m_{v_k},$$ for the vectors $$v_k=v+k\, \langle 1, 0, 1\rangle.$$ The explicit construction is as follows. Since $c_1(v)=H$, by stability it follows that $$H^2(E)=0 \implies h^0(E)-h^1(E)=\chi(v)\leq 0.$$ For each $k\geq -\chi(v)$, consider the Brill-Noether locus $$\m_k=\{E: h^1(E)=k\}\hookrightarrow \m_v$$ and for $E\in \m_k$ construct the universal extension $$0\to H^1(E)\otimes \mathcal O_X\to \widetilde E\to E\to 0.$$ Then, the assignment $$\m_v\ni E\mapsto \widetilde E\in \mathfrak X$$ defines a birational map onto its image. In fact, the fibers of $\pi$ through sheaves $E$ in the Brill-Noether locus $\m_{k}$ are Grassmannians $\mathbb G(k, 2k+\chi(v)).$ 

Now, under the assumptions of the lemma, we have $$\Theta_w^{H^2}=\Theta_{-\lambda_v}^{-\chi(w)}\otimes \Theta_{-\lambda_v-\mu_v}^{s}$$ which shows that $\Theta_w$ is big and nef as well. \qed
\vskip.1in
Finally, the section giving $\Theta[v,w]$ induces a morphism $$\mathsf D:\mathsf V^{\vee}\to \mathsf W.$$ The proposition below suffices to show that $\mathsf D$ is an isomorphism over the generic fiber, thus proving Theorems \ref{genericarb} and 1A.

\begin {proposition}\label{pr3} Over the central fiber, the duality morphism \begin{equation}\label{sdd}\mathsf D_0:H^0(\m_{v_0}, \Theta_{w_0})^{\vee}\to H^0(\m_{w_0},\Theta_{v_0}) \end {equation} is an isomorphism, if $(r,s)\neq (2,2)$ and the dimension inequalities 
\begin{equation}
\label{invpol}
\langle v_0, v_0 \rangle \geq 2(r-1)(r^2+1), \, \,\,\,\langle w_0, w_0 \rangle \geq 2(s-1)(s^2+1)
\end{equation}
 hold. 
 
 When $r=s=2$, the same conclusion is true provided $H^2\geq 8$. \end {proposition} 

\proof Consider a polarization $H_{+}$ suitable with respect to both $v_0$ and $w_0,$ and write $\m_{v_0}^{+}$ and $\m_{w_0}^{+}$ for the moduli spaces of $H_{+}$-semistable sheaves over $\mathcal X_0.$ Theorem \ref{sd} ensures that \begin{equation}\label{d+}\mathsf D^{+}_0:H^0(\m_{v_0}^{+}, \Theta_{w_0})^{\vee}\to H^0(\m_{w_0}^{+}, \Theta_{v_0})\end{equation} is an isomorphism. Assumption \eqref{invpol} is stronger than what is needed to apply the Theorem, provided $(r,s)\neq (2,2)$. When $r=s=2$, Theorem \ref{sd} also holds if $H^2\geq 8$. 

However, \eqref{invpol} is used to apply Corollary \ref{x} and Remark \ref{r4} of the Appendix. In this case, the semistable moduli {\it stacks} do not depend on the choice of polarization away from codimension $2$ loci. Therefore, spaces of sections of theta bundles on the moduli stacks are unaffected by the change of polarization. The translation to the moduli {\it schemes} is straightforward, as lifting sections from the moduli scheme to the moduli stack is an isomorphism, by Proposition $8.4$ in \cite {BL}. To spell out the details, write $$\ms[v]^{ss}\to \Delta,\,\,\, \ms[w]^{ss}\to \Delta$$ for the relative moduli stacks of $\mathcal L$-semistable sheaves, and consider the theta bundles $$\ts_w\to \ms[v]^{ss},\,\,\, \ts_v \to \ms[w]^{ss}.$$ Note the morphisms to the moduli schemes $$p:\ms[v]^{ss}\to \m[v], \,\,\,\,\ p: \ms[w]^{ss}\to \m[w]$$ which match the theta bundles accordingly $$\ts_w=p^{\star} \Theta_w, \,\,\,\, \ts_v=p^{\star} \Theta_v.$$ Lifting \eqref{d+} to the stack, we obtain that $$\mathsf D^{+}_0: H^0(\ms^{+}_{v_0}, \ts_{w_0})^{\vee}\to H^0(\ms^{+}_{w_0}, \ts_{v_0})$$ is an isomorphism. In turn, by the polarization invariance of spaces of theta sections, this shows that $$\mathsf D_0:H^0(\ms^{ss}_{v_0}, \ts_{w_0})^{\vee}\to H^0(\ms^{ss}_{w_0}, \ts_{v_0})$$ is an isomorphism as well. \eqref{sdd} is established descending once again to the moduli scheme. This concludes the proof of the proposition, and of Theorem \ref{genericarb} and 1A along with it.\qed 

\vskip.1in

\section*{\bf APPENDIX: CHANGE OF POLARIZATION FOR MODULI SPACES OF HIGHER RANK SHEAVES OVER $K3$ SURFACES}\vskip.1in
\begin{center} {BY KOTA YOSHIOKA}
\end {center}
\vskip.1in

Let $X$ be a $K3$ surface, and fix a Mukai vector $$v:=(r,\xi, a) \in H^*(X,{\mathbb Z})$$ with
$r>0$. For an ample divisor $H$ on $X$, denote by  $\ms (v), \, \ms_H(v)^{ss},$ and $\ms_H(v)^{\mu-ss}$ the stacks
of sheaves, of Gieseker $H$-semistable sheaves, and of slope $H$-semistable  sheaves respectively -- all of type $v$. 

\begin{lemma}
If $H$ is general with respect to $v$, that is,
$H$ does not lie on a wall with respect to $v$,
then
\begin{equation}
\dim \ms_H(v)^{ss}=
\begin{cases}
\langle v^2 \rangle+1,& \langle v^2 \rangle>0\\
\langle v^2 \rangle+l,& \langle v^2 \rangle=0\\
\langle v^2 \rangle+l^2=-l^2, & \langle v^2 \rangle<0
\end{cases},
\end{equation} 
where
$l=\gcd(r,\xi,a)$.
In particular,
$$\dim \ms_H(v)^{ss} \leq \langle v^2 \rangle+r^2.$$
\end{lemma}

\proof If $\langle v^2 \rangle \geq 0$, then the claims are 
Lemma $3.2$ and $3.3$ in \cite{KY}.
If $\langle v^2 \rangle<0$, then
$\ms_H(v)^{ss}$ consists of $E_0^{\oplus l}$,
where $E_0$ is the unique member of 
$\ms_H(v/l)^{ss}$.
In this case, $\ms_H(v)^{ss}=BGL(l)$,
and $\dim \ms_H(v)^{ss}=-\dim \text{Aut}(E_0^{\oplus l})=-l^2$. 
\qed
\vskip.1in

Let ${\mathcal F}_{H} (v_1,v_2,\dots,v_s)$ be the stack of the Harder-Narashimhan
filtrations
\begin{equation}
0 \subset F_1 \subset F_2 \subset \dots \subset F_s=E, \, \, \, E \in \ms(v)
\end{equation}
such that the quotients
 $F_i/F_{i-1}$, $1 \leq i \leq s$ are semistable with respect to $H$ and 
 \begin{equation}
 v(F_i/F_{i-1})=v_i.
 \label{vi}
 \end{equation}
Then Lemma $5.3$ in \cite{KY} implies 
\begin{equation}
 \dim {\mathcal F}_{H}(v_1,v_2,\dots,v_s)=\sum_{i=1}^s \dim \ms_H(v_i)^{ss}
 +\sum_{i<j} \langle v_i,v_j \rangle.
\end{equation} Note that  $$\text{Hom}(F_i/F_{i-1},F_j/F_{j-1})=0 \, \, \, \text{for}\, \, \,  i<j,$$
as reduced Hilbert polynomials are strictly decreasing in the Harder-Narasimhan filtration.

Let $H_1$ be an ample divisor on $X$ which belongs to
a wall $W$ with respect to $v$ and $H$ an ample
divisor which belongs to an adjacent chamber. 
Then Gieseker $H$-semistable sheaves are $H_1$ slope-semistable $$\ms_{H}(v)^{ss}\hookrightarrow \ms_{H_1}(v)^{\mu\text{-}ss}$$
We shall estimate 
the codimension of 
$$\ms_{H_1}(v)^{\mu\text{-}ss} \setminus \ms_H(v)^{ss}.$$
Specifically, we shall prove
\begin{proposition}\label{prop:H_1}
\begin{equation}
\begin{split}
(\langle v^2 \rangle+1)-
\dim(\ms_{H_1}(v)^{\mu\text{-}ss} \setminus \ms_H(v)^{ss}) 
\geq & \frac{1}{r}\langle v^2 \rangle/2
+r-r^2+1.
\end{split}
\end{equation}
\end{proposition}
As a consequence, we have
\begin{corollary}\label{x}
Assume that
$$\frac{1}{r}\langle v^2 \rangle/2 
+r-r^2+1 \geq 2.$$ Then
$\ms_H(v)^{ss}$ 
is independent on the choice of ample line bundle $H$ (generic or on a wall)
away from codimension $2$. 
\end{corollary}

\proof
Let  $E$ be an $H_1$ slope-semistable sheaf, which is however not  $H$-semistable. Consider its Harder-Narasimhan filtration relative to $H$,
$$0 \subset F_1 \subset F_2 \subset \dots \subset F_s=E.$$ All the subsheaves in the filtration are $H$-destabilizing for $E$. As $E$ is $H_1$ slope-semistable, we must have
equalities of slopes, $$\mu_{H_1} (F_1) = \mu_{H_1} (F_2) = \ldots = \mu_{H_1} (E),$$
or in the notation of \eqref{vi},
\begin{equation}
\label{slopeequal}
\frac{c_1(v_i) \cdot H_1}{\text{rk }v_i}=\frac{c_1(v)\cdot H_1}{\text{rk } v}, \, \,\, \, 1 \leq i \leq s.
\end{equation}
Thus
$$
\ms_{H_1}(v)^{\mu\text{-}ss} \setminus \ms_H(v)^{ss}
=\cup_{v_1,...,v_s}
{\mathcal F}_H(v_1,v_2,\dots,v_s),
$$
where \eqref{slopeequal} is satisfied. We shall estimate 
$\sum_{i<j} \langle v_i,v_j \rangle$.
We set $v_i:=(r_i,\xi_i,a_i)$.
Since
\begin{equation}
\langle (v_i/r_i-v_j/r_j)^2 \rangle=
(\xi_i/r_i-\xi_j/r_j)^2,
\end{equation}
we get
\begin{equation}
\langle v_i,v_j \rangle=
\frac{r_j}{r_i}\langle v_i^2 \rangle/2
+\frac{r_i}{r_j}\langle v_j^2 \rangle/2-
\frac{(r_j \xi_i-r_i \xi_j)^2}{2r_i r_j}. 
\end{equation}
Then we have
\begin{equation}
\begin{split}
\langle v^2 \rangle/2=&
\sum_{i<j} \langle v_i,v_j \rangle+ \sum_i \langle v_i^2 \rangle/2\\
=& \sum_i \frac{r}{r_i}\langle v_i^2 \rangle/2-
\sum_{i<j}\frac{(r_i \xi_j-r_j \xi_i )^2}{2r_i r_j}.
\end{split}
\end{equation}
Hence
\begin{equation}
\sum_i \frac{1}{r_i}\langle v_i^2 \rangle/2=
\frac{1}{r}\langle v^2 \rangle/2+
\frac{1}{r}\sum_{i<j}\frac{(r_i \xi_j-r_j \xi_i )^2}{2r_i r_j}
\end{equation}
and

\begin{equation}
\label{hid}
\begin{split}
\sum_{i<j}\langle v_i,v_j \rangle=&
\sum_i \frac{r-r_i}{r_i}\langle v_i^2 \rangle/2-
\sum_{i<j}\frac{(r_i \xi_j-r_j \xi_i )^2}{2r_i r_j}\\
=&
\sum_i \frac{r-r_i}{r_i}(\langle v_i^2 \rangle/2+r_i^2)
-\sum_i (r-r_i)r_i-
\sum_{i<j}\frac{(r_i \xi_j-r_j \xi_i )^2}{2r_i r_j}\\
\geq & 
\sum_i \frac{1}{r_i}(\langle v_i^2 \rangle/2+r_i^2)
-\sum_i (r-r_i)r_i-
\sum_{i<j}\frac{(r_i \xi_j-r_j \xi_i )^2}{2r_i r_j}\\
=&
\sum_i \frac{1}{r_i}\langle v_i^2 \rangle/2+r-r^2
+\sum_i r_i^2-
\sum_{i<j}\frac{(r_i \xi_j-r_j \xi_i )^2}{2r_i r_j}\\
=&
\frac{1}{r}\langle v^2 \rangle/2+
\frac{1}{r}\sum_{i<j}\frac{(r_i \xi_j-r_j \xi_i )^2}{2r_i r_j}
+r-r^2
+\sum_i r_i^2-
\sum_{i<j}\frac{(r_i \xi_j-r_j \xi_i )^2}{2r_i r_j}\\
\geq &
\frac{1}{r}\langle v^2 \rangle/2
+r-r^2+\sum_i r_i^2,
\end{split}
\end{equation}
where we also used the Hodge index theorem
and Bogomolov's inequality
$$\langle v_i^2 \rangle+2r_i^2 \geq 0.$$
Therefore
\begin{equation}
\nonumber
\begin{split}
(\langle v^2 \rangle+1)-
\dim {\mathcal F}_H(v_1,v_2,\dots,v_s)=&
\sum_{i<j} \langle v_i,v_j \rangle+1-
\sum_i (\dim \ms_H(v_i)^{ss}-\langle v_i^2 \rangle)\\ 
\geq & \frac{1}{r}\langle v^2 \rangle/2
+r-r^2+1,
\end{split}
\end{equation}
which implies the claim.
\qed

\begin {remark} \label{r4} When $c_1(v)$ is primitive, the estimate \eqref{hid} is strict. Indeed, in this case, equality cannot occur in the Hodge index theorem. Therefore, the assumption of Corollary \ref{x} may be relaxed to $$\langle v, v\rangle \geq 2(r-1)(r^2+1).$$
\end {remark}

\end{document}